\documentclass[12pt]{article}
\usepackage{amssymb, amsmath}
\begin{document}
\renewcommand{\theequation}{\arabic{section}.\arabic{equation}}
\newtheorem{thm}{Theorem}[section]
\newtheorem{lem}{Lemma}[section]
\newtheorem{prop}{Proposition}[section]
\newtheorem{cor}{Corollary}[section]
\newcommand{\n}{\nonumber}
\newcommand{\tv}{\tilde{v}}
\newcommand{\tw}{\tilde{\omega}}
\renewcommand{\t}{\theta}
\newcommand{\w}{\omega}
\renewcommand{\b}{\dot{B}^0_{\infty ,1}}
\newcommand{\bn}{\|_{\dot{B}^0_{\infty ,1}}}
\newcommand{\e}{\varepsilon}
\renewcommand{\a}{\alpha}
\newcommand{\La}{\Lambda ^s}
\renewcommand{\l}{\lambda}
\newcommand{\vare}{\varepsilon}
\newcommand{\s}{\sigma}
\renewcommand{\o}{\omega}
\renewcommand{\O}{\Omega}
\newcommand{\bb}{\begin{equation}}
\newcommand{\ee}{\end{equation}}
\newcommand{\bq}{\begin{eqnarray}}
\newcommand{\eq}{\end{eqnarray}}
\newcommand{\bqn}{\begin{eqnarray*}}
\newcommand{\eqn}{\end{eqnarray*}}
\title{Remarks on a Liouville-type theorem for  Beltrami flows}

\author{Dongho Chae$^*$ and Peter Constantin$^\dagger$\\
\ \\
($*$)Department of Mathematics\\
              Chung-Ang University\\
  Seoul 156-756, Republic of Korea\\
 e-mail :dchae@cau.ac.kr\\
 \ \\
 ($\dagger$)Department of Mathematics\\
  Princeton University\\
  Princeton, NJ 08544, USA\\
e-mail: const@math.princeton.edu}
 \date{}
  \maketitle
\begin{abstract}
We present a simple, short and elementary proof that if $v$ is a Beltrami flow with a finite energy in $\mathbb R^3$ then $v=0$. In the case of the Beltrami flows satisfying $v\in L^\infty _{loc} (\Bbb R^3) \cap L^q(\Bbb R^3)$ with $q\in [2, 3)$, or $|v(x)|=O(1/|x|^{1+\vare})$ for some $\vare >0$, we provide a different, simple proof that $v=0$.
  \noindent{\bf
keywords:}   Euler equations, Beltrami flows, Liouville type theorem
\end{abstract}

\section{Introduction}
Ideal homogeneous incompressible inviscid fluid flows are governed by the Euler equations:
 \[
\mathrm{(E)}
 \setcounter{equation}{0}
 \left\{ \aligned
 &\frac{\partial v}{\partial t} +(v\cdot \nabla )v =-\nabla p ,
 \quad (x,t)\in {\mathbb R^3}\times (0, \infty) \\
 &\textrm{div }\, v =0 , \quad (x,t)\in {\mathbb R^3}\times (0,
 \infty)\\
  \endaligned
  \right.
  \]
where $v=(v_1, \cdots, v_n )$, $v_j =v_j (x, t)$, $j=1,\cdots,n$, $n\geq 2$, is the
velocity of the flow, $p=p(x,t)$ is the scalar pressure.
Let $R_j$, $j=1, \cdots, n$, denote the Riesz transforms, given by
 $$R_j(f)(x)= C_n\lim_{\vare\to 0} \int_{\mathbb R^n\setminus B_\vare(x)} \frac{(x_j-y_j)f(y)}{|x-y|^{n+1}} dy,\quad C_n=\frac{\Gamma \left(\frac{n+1}{2} \right)}{\pi ^{\frac{n+1}{2}}}.
 $$
The pressure in the Euler (and Navier-Stokes) equations is given in terms of the velocity up to addition of a harmonic function by
\bb\label{pre-vel}  p= \sum_{j,k=1}^n R_j R_k (v_jv_k).
\ee
This is easily seen by taking the divergence of (E).
In \cite{c1}(see also \cite{bm}) the following result is obtained.
\begin{thm}
If $(v,p)$ satisfies (\ref{pre-vel}) and $|p|+|v|^2 \in L^1(\mathbb R^n)$, then
\bb\label{ortho}
\int_{\mathbb R^n} v_jv_k dx=-\delta_{jk}\int_{R^n} p\,dx.
\ee
\end{thm}
In the next section we present a simple proof of this result using the continuity of the Fourier transform of functions belonging to $L^1 (\mathbb R^n)$. In order to see the implications of the above theorem for Beltrami flows, let us
recall that in the stationary case in $\mathbb R^3$, the first equations of (E) can be rewritten as
\bb
v\times \o=\nabla (p+\frac12 |v|^2),\quad \o=\mathrm{curl}\, v.
\ee
A vector field $v$ is called a Beltrami flow if there exists a function $\lambda=\lambda (x)$ such that
  \bb\label{bel}
  \o =\lambda v.
  \ee
Therefore, if $v$ is a Beltrami flow, then the pair $(v,p)$ is a solution of the stationary Euler equations if
\bb\label{pre}
p+ \frac12 |v|^2 =c, \quad \mbox{$c=$ constant}.
\ee
We call such a solution $(v, p)$ a ``Beltrami solution'' of the stationary Euler equations.
We refer to \cite{ep} for a recent interesting result  regarding  the Beltrami flows.
Recently, Nadirashvili proved a Liouville type property of  Beltrami flows (\cite{n}). He showed that a Beltrami solution $(v,p)$ satisfying
either $v\in L^q(\mathbb R^3)$, $2\leq q\leq 3$, or $|v(x)| =o(1/|x|)$ is necessarily trivial,  $v=0$.
In the case of finite energy Beltrami flows we  have the following immediate consequence of Theorem 1.1:
\begin{thm}
Let $(v,p)$ be a Beltrami solution of the stationary Euler equations with the $\lambda $ given in (\ref{bel}).
If  $v\in L^2(\mathbb R^3)$, then $v=0$. The same conclusion holds, for instance, if  there exists $q\in [\frac65, \infty]$ such that $v\in L^q(\mathbb R^3)$ and $\lambda \in L^{\frac{6q}{5q-6}} (\mathbb R^3)$ (if $v\in L^{\frac65} (\mathbb R^3)$, then we require $\lambda \in L^\infty (\mathbb R^3)$ ).
\end{thm}
We have also the following result for the cases considered in the paper \cite{n}, for which we present a different, simple proof.
\begin{thm}
Let $v\in L^\infty _{loc} (\Bbb R^3)$ be a Beltrami solution of the stationary Euler equations satisfying
either $v\in L^q (\Bbb R^3)$ for some $q\in [2, 3)$, or that there exists $\vare >0$ such that $ |v(x)|= O(1/|x|^{1+\vare})$ as $|x|\to \infty$.
 Then, $v=0$.
 \end{thm}
\section{Proof of the Theorems}
 \setcounter{equation}{0}
 We use the notation for the Fourier transform  of $f(x)$
 $$\mathcal{F}(f)(\xi)=\widehat{f}(\xi)=\frac{1}{(2\pi)^{\frac{n}{2}}} \int_{\mathbb R^n} f(x)e^{-ix\cdot \xi } dx,
 $$
 whenever the right hand side is defined.
 In terms of the Fourier transform the Riesz transform is defined as
 $$ \widehat{R_j (f)}(\xi)=\frac{i \xi_j}{|\xi|}, \quad i=\sqrt{-1}.
 $$
\noindent{\bf Proof of Theorem 1.1 } Without loss of generality we may restrict ourselves  to the stationary case, $v(x,t)=v(x), p(x,t)=p(x)$. By the Fourier transform one has
\bb\label{th1} \hat{p}(\xi)=-\sum_{j,k=1}^n \frac{\xi_j\xi_k}{|\xi|^2} \widehat{v_j v_k }(\xi).
\ee
We note that $\widehat{p}(\xi)$ and $\widehat{v_j v_k}(\xi)$, $j, k =1,  \cdots, n$, are continuous at $\xi=0$ from the hypothesis, $|p|+|v|^2\in L^1(\mathbb R^n)$.
  Let $w$ be a given constant vector with $|w|=1$. We put $\xi= \rho w $ in (\ref{th1}), and pass $\rho \to 0$ to obtain
\bb\label{th2}
\int_{\mathbb R^n} p\, dx = -\int_{\mathbb R^n} (v\cdot w )^2 dx.
\ee
If we plug $w=\mathbf{e}^j$ in (\ref{th2}), where $\mathbf{e}^j$ is the canonical basis of $\mathbb R^n$ with its components given by $(\mathbf{e}^j)_k =\delta_{jk}$, then we have
$$\int_{\mathbb R^n} p \, dx= -\int_{\Bbb R^n} v_j ^2 dx\quad \forall j=1, \cdots, n.
$$
On the other hand, for $j\neq k$, if we  put $w=\frac{\mathbf{e}^j+\mathbf{e}^k}{\sqrt{2}}$ in (\ref{th2}), we obtain
$\int_{\mathbb R^n} v_j v_k dx=0$.
 $\square$\\
 \ \\
 \noindent{\bf Proof of Theorem 1.2 } Since $(v,p)$ is a Beltrami solution of stationary Euler equations, we have $p-c=\frac12 |v|^2:=\tilde{p}$ for some constant $c$.
 In the case $v\in L^2 (\mathbb R^3)$ we find that $|v|^2+|\tilde{p}|\in L^1(\mathbb R^3)$, and by Theorem 1.1 we obtain
 $$\int_{\mathbb R^3} \tilde{p} \, dx= -\frac{1}{3}\int_{\mathbb R^3} |v|^2 dx=-\frac12\int_{\mathbb R^3} |v|^2 dx,
 $$
 which implies that $v=0$, and  $\tilde{p}=\sum_{j,k=1}^n R_j R_k (v_jv_k)=0$.  On the other hand, if $v\in L^{q} (\mathbb R^3)$ and $\lambda \in L^{\frac{6q}{5q-6}} (\mathbb R^3)$ with $\frac65<q\leq \infty$, or $v\in L^{\frac65} (\mathbb R^3)$ and $\lambda \in L^\infty (\mathbb R^3)$, then
 we estimate
 \bqn
 \|v\|_{L^2}&\leq& C \|\nabla v\|_{L^{\frac{6}{5}}}\leq C \| \o\|_{L^{\frac{6}{5}}}=C \| \lambda v \|_{L^{\frac{6}{5}}}\\
 &\leq& C\|\lambda\|_{L^{\frac{6q}{5q-6}}} \|v\|_{L^q} <\infty,
 \eqn
 and we reduce to the above case of $v\in L^2(\Bbb R^3)$. $\square$\\
 \ \\
 \noindent{\bf Proof of Theorem 1.3 }
 We first observe that our hypothesis implies that
 \bb\label{th13a} \int_{\Bbb R^3} |v|^2 |x| ^{\lambda -2} dx <\infty.
 \ee
 for some $\lambda \in (1, 2)$.
 Indeed, the case $ |v(x)|= O(1/|x|^{1+\vare})$ as $|x|\to \infty$ is obvious, while in the case $v\in L^q (\Bbb R^3)$ for some $q\in [2, 3)$, we have the following estimate,
 \bqn
 \int_{\{|x|\geq 1\}} |v|^2|x|^{\lambda -2} dx\leq C \|v\|_{L^q}^2 \left(\int_1 ^\infty  r^{2-\frac{q(\lambda -2)}{q-2}} dr\right)^{\frac{q-2}{q}}<\infty
 \eqn
 for $\lambda$ with  $1< \lambda < \frac{6}{q}-1$.
 Let us introduce a standard radial cut-off function $\sigma\in
C_0 ^\infty(\Bbb R^N)$ such that
 \bb\label{16}
   \sigma(|x|)=\left\{ \aligned
                  &1 \quad\mbox{if $|x|<1$},\\
                     &0 \quad\mbox{if $|x|>2$},
                      \endaligned \right.
 \ee
and $0\leq \sigma  (x)\leq 1$ for $1<|x|<2$.  Then, for each $R
>0$, we define
 $
\s \left(\frac{|x|}{R}\right):=\s_R (|x|)\in C_0 ^\infty (\Bbb R^N).
$
A Beltrami solution $(v,p)$ with $p=-\frac12 |v|^2 +C$ satisfies
\bb\label{th13b}
\sum_{j,k=1}^3 \int_{\Bbb R^3} v_jv_k \partial_j\partial_k \varphi dx=\frac12 \int_{\Bbb R^3} |v|^2 \Delta \varphi dx \quad \forall \varphi \in C_0 ^2 (\Bbb R^3).
\ee
We choose our test function $ \varphi (x)=\varphi_{\delta, R}(x)= (|x|^{2\lambda } +\delta)^{\frac12} \s_R$ for $\delta, R>0$ in (\ref{th13b}), which is an approximation of $\varphi =|x|^\lambda$, and passing
first $\delta \to 0$, and then $R\to \infty$, using the dominated convergence theorem, taking (\ref{th13a}) into account, we obtain easily that
\bb\label{th13c}
(\lambda -1)\int_{\Bbb R^3} |v|^2 |x|^{\lambda -2} dx=2(\lambda -2) \int_{\Bbb R^3} (v\cdot x)^2 |x|^{\lambda -4} dx.
\ee
The fact that $\lambda \in (1,2)$ in (\ref{th13c}) implies $v=0$. $\square$

 $$\mbox{\bf Acknowledgements}$$
DC was partially supported by NRF grants 2006-0093854 and  2009-0083521.
PC was partially supported by NSF DMS grants 1209394 and 1265132.

\end{document}